\RequirePackage{ifpdf}
\ifpdf % We are running pdfTeX in pdf mode
\documentclass[pdftex]{sigma}
\else
\documentclass{sigma}
\fi

\def\g{{\cal{G}}}
\def\k{{\mathbf{k}}}
\def\r{\mathbb{R}}

\begin{document}

\renewcommand{\PaperNumber}{088}

\FirstPageHeading

\ShortArticleName{Geodesically Complete Lorentzian Metrics on Some Homogeneous 3 Manifolds}

\ArticleName{Geodesically Complete Lorentzian Metrics\\ on Some Homogeneous 3 Manifolds}

\Author{Shirley BROMBERG~$^\dag$ and Alberto MEDINA~$^\ddag$}

\AuthorNameForHeading{S. Bromberg and A. Medina}

\Address{$^\dag$~Departameto de Matem\'aticas, UAM-Iztapalapa,
M\'exico}

 \EmailD{\href{mailto:stbs@xanum.uam.mx}{stbs@xanum.uam.mx}}

\Address{$^\ddag$~D\'epartement des Math\'ematiques, Universit\'e de
Montpellier II, UMR, CNRS,\\
\hphantom{$^\ddag$}~5149, Montpellier, France}
\EmailD{\href{mailto:medina@math.univ-montp2.fr}{medina@math.univ-montp2.fr}}

\ArticleDates{Received June 24, 2008, in f\/inal form December 10,
2008; Published online December 18, 2008}

\Abstract{In this work it is shown that a necessary condition for
the completeness of the geodesics of left invariant
pseudo-Riemannian metrics on Lie groups is also suf\/f\/icient in the
case of 3-dimensional unimodular Lie groups, and not suf\/f\/icient for
3-dimensional non unimodular Lie groups. As a consequence it is
possible to identify, amongst the compact locally homogeneous
Lorentzian 3-manifolds with non compact (local) isotropy group,
those that are geodesically complete.}


\Classification{53C22; 53C50; 57M50; 22E30}

\vspace{-2mm}

\section{Introduction}

The aim of this paper is to characterize the geodesically complete
left invariant Lorentzian metrics on unimodular 3-dimensional Lie
groups. The main result (Theorem~\ref{theorem2}) improve and
complete Proposition~3.2 given in \cite{BM1997}. It is well known
that an unimodular Lie group of dimension~3 contains an inf\/inity of
cocompact discret subgroups. Consequently our results provide a~rich
class of compact and complete (or incomplete) Lorentzian
3-manifolds. In particular, we obtain the geodesically complete and
compact locally homogeneous Lorentz 3-manifolds with non compact
(local) isotropy group as described in~\cite{DZ}. Moreover, our
results give a correct description of the complete (or incomplete)
Lorentzian metrics on $M=\Gamma\backslash {\mathrm{PSL}}(2,\r),$
where $\Gamma$ is a Fuchsian group, obtained from a left invariant
Lorentzian metric on $\mathrm{PSL}(2,\r).$

To characterize the complete left invariant pseudo-Riemannian
metrics which are not def\/inite on arbitrary Lie groups is a
dif\/f\/icult problem. On one side of the spectrum we f\/ind groups where
all such metrics are complete, such as the Abelian groups, the
compact groups, the 2-nilpotent groups \cite{G1994} or the group of
isometries of the Euclidean plane, and on the other side, groups
such as the group of the af\/f\/ine transformations of the line in which
no Lorentzian metric is complete.

The structure of this work is as follows:
In Section~\ref{section1} the def\/initions and some general results are given,
particularly the algebraic counterparts of the geometric properties
of left invariant pseudo-Riemannian metrics def\/ined on Lie groups.
In Section~\ref{section2} we characterize the left invariant metrics on
3-dimensional unimodular Lie groups that are geodesically complete.
The following result is proved:

\medskip

 %{\bf{Theorem 2}}
 {\it A left invariant
metric on an unimodular $3$-dimensional Lie group is geodesically
complete if and only if all geodesics that run on $1$-parameter
subgroups are complete. }

\medskip

 In Section~\ref{section3} we give an example that shows that
Theorem~\ref{theorem2} is not valid without the hypothesis of
unimodularity. We also show that previous results in the literature~\cite{GL} are not accurate.

\section{Preliminaries and general results}\label{section1}

Let $\g={\mathrm{T}}_\epsilon(G)$ be the Lie algebra of a Lie group
$G,$ where $\epsilon$ stands for the identity of $G.$ Every left
invariant pseudo-Riemannian metric $q$ on $G$ induces a left
invariant connection $\nabla$  on $G,$ the Levi-Civita connection; a
non degenerate quadratic form $\langle \cdot ,\cdot \rangle$ on $\g$, the polar
form of $q_\epsilon$, the restriction of $q$ to $\g $; and a product
$x\cdot y$ on $\g$ def\/ined by
\[
(x\cdot y)^+=\nabla_{x^+}{y^+},
\] where $x^+$ is the left invariant
vector f\/ield on $G$ determined by $x\in\g.$ These objects are
related by the Koszul formula:
\[
\langle x\cdot
y,z\rangle=\tfrac{1}{2}\left(\langle[x,y],z\rangle-\langle[y,z],x\rangle+\langle[z,x],y\rangle\right).
\]

If $c(t)$ is a curve in $G$ and
\[
x(t)=\left(L_{c(t)}\right)^{-1}_{*,c(t)}c'(t),\] where $L_\sigma$
denotes the left multiplication by $\sigma$ in $G,$ it follows that
$c(t)$ is a geodesic if and only if $x(t)$ is a solution of the
equation
\begin{equation}\label{equation1}
\dot{x}=-x\cdot x .\end{equation}

That $\nabla$ is geodesically complete is equivalent to the
completeness of the vector f\/ield given by equation
(\ref{equation1}). Clearly  the {\it energy}, that is
$q_\epsilon(x)=\langle x,x\rangle,$ is a f\/irst integral of equation~(\ref{equation1}).

Notice that if $G$ is Abelian, then $x\cdot y\equiv 0$ for every
left invariant pseudo-Riemannian metric. Hence, every left invariant
pseudo-Riemannian metric on an Abelian Lie group is complete.

If $\phi :\g\to\g^{*}$ is the canonical linear isomorphism
associated to $\langle\cdot,\cdot\rangle$ then the equation
(\ref{equation1}) can be written, in terms of the co-adjoint
representation as
\begin{equation}\label{equation2}
\dot{\xi}=-{\mathrm{ad}}^{*}_{\phi^{-1}\xi}\xi
\end{equation}
with $\xi\in\g^{*}.$  Equation (\ref{equation2}) is called the Euler
equation associated to the pseudo-Riemannian metric, and the
corresponding vector f\/ield is called the Euler vector f\/ield.
Consequently the solutions of equation (\ref{equation2}) are curves
in the orbits of the co-adjoint representation of $G.$ This implies
that (\ref{equation2}) is complete if $G$ is compact. The energy on
$\g^{*} $ is given by $q_\epsilon(\phi^{-1}\xi).$

Suppose that $(G,\k)$ is orthogonal or quadratic, i.e.\ $\k$ is a
bi-invariant pseudo-Riemannian metric. Then $x\cdot x=0$ and
equation (\ref{equation1}) is obviously complete and the geodesics
through $\epsilon$ are the 1-parameter subgroups of $G.$ Note that
every left invariant pseudo-Riemannian metric~$\langle\cdot,\cdot\rangle$
on a orthogonal Lie group $(G,\k)$ is given by the formula $\langle
x,y\rangle=\k(u(x),y)$ where $u$ is a~$\k$-symmetric linear
isomorphism of $\g.$ Equation (\ref{equation1}) is given in this
case by
\[u(\dot{x})=[x,u(x)]\]
which is equivalent by a linear change of coordinates to the Lax
pair
\begin{equation}\label{equation3}\dot{x}=[x,u^{-1}x].\end{equation}
Equation (\ref{equation3}) has two
quadratic f\/irst integrals $\k(x,x)$, $\k(x,u^{-1}x).$

The vector f\/ields induced by a left invariant pseudo-Riemannian
metric are special cases of {\it homogeneous quadra\-tic vector
fields.}

\begin{definition}
A {\it homogeneous quadratic vector field} on the $n$-dimensional
af\/f\/ine space~$\r^n$ is a~vector f\/ield $F:\r^n\to\r^n,$
$F=(F_1,\ldots, F_n),$ where $F_i$ is a homogeneous polynomial of
degree~2, $i=1,\ldots, n.$ If the solutions of the dif\/ferential
equation $\dot{X}=F(X)$ are def\/ined for all $t\in\r,$ the vector
f\/ield $F$ is called {\it complete}.
\end{definition}

The following def\/inition provides the simplest example of complete
homogeneous quadratic vector f\/ields:

\begin{definition}
An {\it affine-quadratic} (or linearizable) vector f\/ield is a
homogeneous qua\-dratic vector f\/ield such that $\dot{X}=F(X)$ is
equivalent by a linear change of coordinates to
\[
\dot{U}=A(V)U+B(V), \qquad \dot{V}= 0,
\]
where $A$ is linear and the coordinate functions of $B$  are
homogeneous polynomials of degree~2.
\end{definition}

 It is easily checked that every
af\/f\/ine-quadratic vector f\/ield is complete.

\begin{definition} Let $F$ be a homogeneous quadratic vector f\/ield. We def\/ine
\[
{\cal{I}}:=\{X\ne 0\, ;\, X, F(X)\,\, {\mbox{{\rm{are linearly
dependent}}}}\}.
\] If $X\in{\cal{I}}$ and $F(X)=0,$ $X$ is called a
{\it zero} of $F.$ If $X\in{\cal{I}}$ and $F(X)\ne 0$ then $X$  is
called a~relaxed {\it idempotent} of $F.$
\end{definition}

  In fact a strict idempotent is a non trivial
solution of $F(X)=X.$ If $X\in{\cal{I}}$ and $F(X)\ne 0,$ then there
exists $\lambda\ne 0$ such that $F(\lambda X)=\lambda X,$ i.e.\ for
every relaxed idempotent there is a~strict idempotent in the same
direction (the line through $X$ is invariant by $F$). In what
follows we make no distinction between relaxed or strict
idempotents.

\begin{remark}\label{remark1}A simple calculation shows that any idempotent of
$F$ gives an incomplete solution for $\dot{X}=F(X).$
\end{remark}

  In fact the solution of $\dot{X}=F(X)$ with initial
condition an idempotent $X_0$ is given by $t\mapsto \alpha (t)X_0$
where $\alpha$ satisf\/ies $\dot\alpha=\alpha^2.$ Moreover if $X_0\in
{\cal{I}}$ for the Levi-Civita product, the geodesic through~$\epsilon$ with velocity $X_0$ is given by $t\mapsto
\exp \beta(t)X_0,$ where $\beta$ is such that $\beta(0)=0,$
$\dot{\beta}(t)=\alpha (t)$ and $\dot\alpha=\alpha^2.$ A result by
Kaplan and Yorke (see~\cite{KY}) states that for every quadratic
homogeneous vector f\/ield $F$, ${\cal{I}}\ne \varnothing$, hence

\begin{proposition}\label{proposition1}
Any left invariant pseudo-Riemannian metric on $G$ has a geodesic
that runs on a $1$-parameter subgroup. Furthermore if the
pseudo-Riemannian metric is complete then the corresponding
Levi-Civita product has a zero. Hence for such metrics, there is a
$1$-parameter subgroup of $G$ which is a geodesic.
\end{proposition}

\begin{lemma}\label{lm:idemp} Let $X$ be an  idempotent of a quadratic
vector field $F$ and $Q$ a quadratic form which is a first integral
of $\dot{X}=F(X).$ Then $Q(X)=0.$
\end{lemma}

\begin{proof}
Let $P$ be the bilinear polar form associated to $Q.$ Then for any
curve $t\mapsto X(t)$
\[ \frac{d}{dt}Q(X)=
2P\left(X,\frac{d}{dt}X\right).
\]
If $t\mapsto X(t)$ is a solution of $\dot{X}=F(X)$ and $Q$ is a
f\/irst integral, then
\[
0=\frac{d}{dt}Q(X)= 2P\left(X,\frac{d}{dt}X\right)= 2P(X,F(X))=-2Q(X)
\]
as was asserted.
\end{proof}

 The following theorem is a general result of
completeness for homogenous quadratic vector f\/ields in dimension~2.

\begin{theorem}[\cite{BM2005}]\label{theorem:QTDS}
Let $F$ be a homogeneous quadratic vector field on the plane. The
following assertions are equivalent:
\begin{enumerate}\itemsep=0pt
\item[$(i)$] $\dot{X}=F(X)$ is complete;

\item[$(ii)$] $F$ is quadratic-affine or the differential equation $\dot{X}=F(X)$
is equivalent, by a linear change of coordinates, to a system
\[
\dot{x}=y(ax+by), \qquad \dot{y}= y(cx+dy),
\]
where $(a+d)^2-4(ad-bc)<0.$
\end{enumerate}
\end{theorem}

\section[Complete left invariant pseudo Riemannian metrics defined on 3-dimensional
unimodular Lie groups]{Complete left invariant pseudo Riemannian metrics def\/ined\\ on 3-dimensional
unimodular Lie groups}\label{section2}

The Lie bracket of a 3-dimensional Lie algebra can be described
using the vector product $\times$ corresponding to the positive
def\/inite scalar product $\langle\cdot ,\cdot\rangle$ and to an orientation
as follows~\cite{M}:
\[
[X,Y]=L(X\times Y),
\]
where $L\in\rm{End}\,(\g).$

If the group is unimodular, i.e.\ $\mathrm{tr}({\rm{ad}}_X )=0$
for all  $X\in\g,$ then $L$ is  $\langle\cdot,\cdot\rangle$-symmetric. Hence
$\g$ has a $\langle\cdot,\cdot\rangle$-orthonormal basis of eigenvectors of
$L,$ ${\cal{B}}=(E_1,E_2,E_3),$ that is
\[ [E_2,E_3]=\alpha_1\,E_1,\qquad [E_3,E_1]=\alpha_2\, E_2,\qquad
[E_1,E_2]=\alpha_3\,E_3. \]

The signature of the eigenvalues characterizes the Lie algebra as is
shown in the following diagram

\centerline{\begin{tabular}{|c||c|c|c|}\hline
\rule[-2mm]{0cm}{6mm} &signature of $\alpha_1$, $\alpha_2$, $\alpha_3$   &
corresponding Lie group & type \\ \hline
\rule[-2mm]{0cm}{5mm}(a)&$0,0,0$ & $\r\oplus\r\oplus\r$  & Abelian
\\
\hline \rule[-2mm]{0cm}{5mm}(b)&$+$, $0$, $0$    & Heisenberg           &
nilpotent  \\
\hline
\rule[-2mm]{0cm}{6mm}(c)&$+$, $+$, $+$    & $\mathrm{O}(3)$         & compact  \\
\hline \rule[-2mm]{0cm}{5mm}(d)&$+$, $+$, $0$    & $\rm{E}(2)$
& solvable  \\
\hline \rule[-2mm]{0cm}{5mm}(e)&$+$, $-$, $0$    & $\mathrm{E}(1,1)$ &
solvable
\\
\hline \rule[-2mm]{0cm}{5mm}(f)&$+$, $+$, $-$    & $\mathrm{Sl}(2,\r)$ &
semisimple  \\ \hline
\end{tabular}}

\begin{theorem}\label{theorem2}
The following assertions are equivalent for left invariant
pseudo-Riemannian metrics defined on $3$-dimensional unimodular Lie
groups
\begin{itemize}\itemsep=0pt
\item[$(i)$] the Euler equation is complete;

\item[$(ii)$] all geodesics that run on $1$-parameter subgroups are
complete;

\item[$(iii)$] the Euler field has no idempotents.
\end{itemize}

\end{theorem}

  Notice that, by Remark~\ref{remark1},
$(ii)\Rightarrow (iii)$ since an idempotent of the Euler equation
produces an incomplete geodesic that runs on a 1-parameter subgroup.
Hence the only implication that requires proof is
$(iii) \Rightarrow (i)$.

  In the f\/irst four cases, all pseudo-Riemannian metrics
are complete, the f\/irst three being widely known:
\begin{itemize}\itemsep=0pt
\item[(a)] the Euler equation for a left invariant pseudo-Riemannian metric on an
Abelian Lie group is complete;

\item[(b)] the Euler equation for a left invariant pseudo-Riemannian metric on
Heisenberg group is linearizable;

\item[(c)] the Euler equation for a
left invariant pseudo-Riemannian metric on a compact Lie group has a
def\/inite quadratic f\/irst integral: the Killing form.
\end{itemize}

The proof of Theorem \ref{theorem2} consists on considerations of
case by case situations.

\begin{proposition}\label{proposition2}
Every left invariant pseudo-Riemannian metric on ${\mathrm{E}}(2)$
is complete.
\end{proposition}

\begin{proof} In a suitable basis, the Lie bracket is given by the
relations: \begin{equation} [E_2,E_3]=0,\qquad [E_3,E_1]=E_2,\qquad
[E_1,E_2]=E_3 \label{equation3.2}
\end{equation} and the Euler equation is
\begin{gather*}
\dot{\xi}_1= x_3\,\xi_2 - x_2\xi_3, \qquad \dot{\xi}_2 =  x_1\xi_3, \qquad \dot{\xi}_3=- x_1\xi_2.
\end{gather*}

Since the semi-def\/inite quadratic  $\xi_2^2 +\xi_3^2$ is a f\/irst
integral of the Euler vector f\/ield, by Lemma~\ref{lm:idemp}  the
Euler vector f\/ield has no idempotents.

We split the proof of the completeness of Euler vector f\/ield in two
cases:

  $i)$ If the pseudo-Riemannian metric $\langle\cdot ,\cdot\rangle$
is degenerate on the vector space spanned by~$E_2$,~$E_3$, then $x_1$
does not depend on $\xi_1$ and the two last equations of the Euler
equation lie in the $\xi_2\xi_3$ plane. This planar system of
equations has a positive def\/inite quadratic f\/irst integral, hence it
is complete. The f\/irst equation is of the form \[\dot{\xi}_1=
f(t)\xi_1 + g(t)\]and its solution is def\/ined for all values of $t$
since $f$ and $g$ are so.

  $ii)$ If the pseudo-Riemannian metric $\langle\cdot ,\cdot\rangle$
is non degenerate on  the vector space spanned by~$E_2$,~$E_3$, then
there is a basis of $\mathrm{E}(2)$ with $E_1\perp_{\langle\cdot
,\cdot\rangle} {\mathrm{Vect}}\{E_2,E_3\}$ that satisf\/ies the same
rela\-tions for the bracket stated in equation~(\ref{equation3.2}).
Then the f\/irst integral $\xi_2^2 +\xi_3^2$ and the energy, that in
this basis has the form  $\xi_1^2+a\xi_2^2 +b\xi_2\xi_3+c\xi_3^2,$
can be combined to get a positive def\/inite f\/irst integral. Hence the
solutions of the Euler equation are complete.
\end{proof}

  In \cite{G1996} another proof of Proposition~\ref{proposition2} and an alternative statement and proof of
Proposition~\ref{pr:E(1,1)} are provided.

\begin{proposition} \label{pr:E(1,1)}
If the Euler product of a pseudo-Riemannian metric defined on
$E(1,1)$ has no idempotents then the Euler equation has a definite
quadratic first integral, and hence such a~metric is complete.
\end{proposition}

\begin{proof}
In a suitable basis, the Lie bracket is given by the relations
\begin{equation}\label{eq:brE11}
[E_1,E_2]=E_2,\qquad [E_1,E_3]=-E_3,\qquad
[E_2,E_3]=0 \end{equation} and the Euler equation is
\begin{gather*}
\dot{\xi}_1 = - x_2\,\xi_2 + x_3\xi_3,\qquad
\dot{\xi}_2 =  x_1\xi_2,\qquad
\dot{\xi}_3 = - x_1\xi_3 .
\end{gather*}

Besides the energy:
$\langle\Phi^{-1}\xi,\Phi^{-1}\xi\rangle=x_1\,\xi_1+x_2\,\xi_2+x_3\,\xi_3,$
the Euler vector f\/ield has $\xi_2\xi_3$ as a f\/irst integral. If the
Euler vector f\/ield has no idempotents, the pseudo-Riemannian metric
is non degenerate on $[\g,\g]={\mathrm{Vect}}\{E_2,E_3\}.$ In fact,
for a pseudo-Riemannian metric which is degenerate on $[\g,\g]$,
$x_1$ does not depend on $\xi_1$ and any $\xi_2$ such that $x_1=1$
for $\xi_3=0$ gives an idempotent of the Euler vector f\/ield.

Hence, if the Euler vector f\/ield has no idempotents, there exists
$e_1\perp_{\langle\cdot,\cdot\rangle}[\g,\g]$ such that the bracket
relations given by equation (\ref{eq:brE11}) are still satisf\/ied,
and the energy in these coordinates is given by:
\[
\lambda\xi_1^2+\mu\xi_2^2+2c\xi_2\xi_3+\nu\xi_3^2.
\]
Taking $\xi_2=0,$ the equation for idempotents becomes
\begin{gather*}
\xi_1 =   x_3\xi_3= \nu\xi_3^2,\qquad
\xi_3 = - x_1\xi_3=-\lambda\xi_1\xi_3 .
\end{gather*}
which has a solution if and only if  $\lambda\mu<0.$ Hence, if the
Euler equation has no idempotents $\lambda\mu>0.$ Analogously
$\lambda\nu>0.$ Since $\xi_2\xi_3$ is a f\/irst integral
\[
\lambda\xi_1^2+\mu\xi_2^2+\nu\xi_3^2
\]
is a f\/irst integral which is a def\/inite quadratic form.
\end{proof}

\begin{proposition}\label{proposition4}
If the Lax pair vector field  for a pseudo-Riemannian metric on
${\mathrm{Sl}}(2,\r)$ has no idempotents then it is complete.
\end{proposition}

\vspace{.1in} In order to deal with ${\mathrm{Sl}}(2,\r)$ some
previous considerations are needed. The Killing form~$\k$ of~${\mathit{sl}}(2,\r)$ is a quadratic form of index $(2,1),$ and
induces a bi-invariant Lorentz metric on~${\mathrm{Sl}}(2,\!\r)$. Thus
the completeness of a left invariant pseudo-Riemannian metric is
equivalent to the completeness of the Lax pair. We express this
equation in two adapted basis (all non stated products are obtained
either by symmetry or by antisymmetry):

\begin{itemize}\itemsep=0pt
\item[(A)] given a $\k$-orthonormal basis, i.e.\ an ordered basis
$(E_1,E_2,E_3),$ such that
\[
 1=k(E_1,E_1)= k(E_2,E_2)= -k(E_3,E_3),
\]
the Lie bracket is given by the table
\[
\begin{array}{c|c|c|c}
\rule[-2mm]{0cm}{4mm}[\, ,\,]  & E_1  & E_2   & E_3 \\ \hline
\rule[-2mm]{0cm}{6mm} E_1
&  0   & -\kappa\, E_3  & -\kappa\, E_2 \\
\hline \rule[-2mm]{0cm}{6mm} E_2
&\phantom{-\kappa\, E_3}  &      0               & \phantom{-1}\kappa\, E_1\\
\hline \rule[-2mm]{0cm}{6mm} E_3      &      & &         0
\end{array}
\]
and  matrix of ${\mathrm{ad}}_X$ is
\[
{\rm ad}_X=\kappa\left(\begin{array}{ccc}
     0        & -\, x_3   &  x_2 \\
 x_3  &   0            & -x_1\\
x_2  & -x_1    &     0
\end{array}\right) ;
\]

\item[(B)] given a  $\k$-hyperbolic
basis, i.e.\ an ordered basis $(E_1,E_2,E_3),$ such that
\[
k(E_1,E_1)=k(E_2,E_3)=1,
\]
the Lie bracket is given by
\[
\begin{array}{c|c|c|c}
\rule[-2mm]{0cm}{4mm}[\, ,\,]  & E_1  & E_2   & E_3 \\ \hline
\rule[-2mm]{0cm}{6mm} E_1
&  0   & \kappa\, E_2  & -\kappa\, E_3 \\
\hline \rule[-2mm]{0cm}{6mm} E_2
&\phantom{\epsilon_3\kappa\, E_3}  &      0   & \phantom{-1}\kappa\, E_1\\
\hline \rule[-2mm]{0cm}{6mm} E_3      &      & &         0
\end{array}
\]
and ${\mathrm{ad}}_X$ has the matrix
\[
{\rm ad}_X=\kappa\left(\begin{array}{ccc}
 0    &  -x_3  & x_2    \\
-x_2  &   x_1  & 0\\
x_3   &   0    & -x_1
\end{array}\right),
\]
where $\kappa = k([E_1,E_2],E_3).$ Since it is possible to take care
of a constant by a linear change of coordinates, it will be omitted
in what follows.
\end{itemize}

As was said before, the proof will consist in dealing with each case
of the Jordan form of~$u$.

\begin{lemma}
If $u$ has exactly one eigendirection then the Lax pair vector field
has always idempotents.
\end{lemma}

\begin{proof} Let $E$ be an eigenvector. There are two
possibilities: either
 $\k(E,E)\ne 0$ or $\k(E,E)=0.$ Assume with no loss of generality that
 the eigenvalue is~1.

Suppose f\/irst that $E$ is an eigenvector such that $\k(E,E)\ne 0.$
Then $\k(E,E)>0.$ In fact, under the hypothesis, the subspace
$E^\perp$ is a 2-dimensional non degenerate subspace  which is $u$~invariant. If $\k(E,E)<0,$ $\k$ is positive def\/inite on $E^\perp.$
Thus $u$ is diagonalizable on $E^\perp,$ which contradicts the
assumption.

Let $E_1$ be an eingenvector such that $\k(E_1,E_1)=1,$ and
$E_2,E_3\in E^\perp$ such that $\k(E_2,E_2)=-\k(E_3,E_3)=1,$
$\k(E_2,E_3)=0.$ The matrix of $u$ and of $u^{-1}$ are:
\[
\left(\begin{array}{ccc}
 1&0&0\\
0&\alpha &\beta\\
0&-\beta &\delta
\end{array}\right),\qquad
\frac{1}{\Delta}\left( \begin{array}{ccc}
\Delta&0&0\\
0&\delta&-\beta\\
0&\beta &\alpha
\end{array}\right),
\]
where $\Delta=\beta^2+\alpha\delta.$ Notice that in this case the
characteristic polynomial of $u,$ $p_u(x)=(x-1)$ $(x-(\alpha +\delta
)x+\Delta),$ has one simple real root. Hence $\Delta>0$ and
$(\alpha-\delta)^2-4\beta^2<0.$ The Lax pair is given by:
\begin{gather*}
\dot{x}_1= \frac{1}{\Delta}(\beta x_2^2+(\alpha-\delta)x_2x_3+\beta x_3^2),\\
\dot{x}_2= \frac{x_1}{\Delta}(-\beta x_2+(\Delta-\alpha)x_3),\\
\dot{x}_3= \frac{x_1}{\Delta}((\Delta-\delta) x_2+\beta x_3).
\end{gather*}

If $(x_1,x_2,x_3)$ is an idempotent of the equation  then $x_1\ne 0$
and $\Delta/x_1$ is an eigenvalue for
\[
\left(\begin{array}{cc}
-\beta & \Delta-\alpha\\
\Delta-\delta&\beta
\end{array}
\right).
\]
Since the characteristic polynomial of this matrix is
$x^2-(\beta^2+(\Delta-\alpha)(\Delta-\delta))$, then the planar
system has idempotents if and only if
\[
0<\beta^2+(\Delta-\alpha)(\Delta-\delta)=\Delta(\beta^2+\alpha\delta
+1-(\alpha+\delta))
\]
or if and only if $ \beta^2+\alpha\delta +1-(\alpha+\delta)>0,$
because $\Delta>0.$ This inequality holds since
\[
\beta^2+\alpha\delta > \tfrac{1}{4}(\alpha+\delta)^2\ge
(\alpha+\delta)-1.
\]

Hence the planar system has two eigendirections corresponding to two
eigenvalues of dif\/ferent sign. Finally to produce an idempotent for
the Lax pair, it is necessary to solve the equation
\[
\frac{\Delta}{x_1}= \frac{\Delta^2}{\beta x_2^2+(\alpha-\delta)x_2x_3+\beta x_3^2},
\]
where $\Delta/x_1$ is an eigenvalue and $(x_2,x_3)$ a corresponding
eigenvector. Since the quadratic form $\beta
x_2^2+(\alpha-\delta)x_2x_3+\beta x_3^2$ is def\/inite, this equation
has always a solution.

 Now suppose that $\k(E,E)=0.$ Let $E_3$ be an
eigenvector. By hypothesis, $\k(E_3,E_3)=0.$ Let~$E_2$ be such that
$\k(E_3,E_2)=1$ and $\k(E_2,E_2)=0.$ Choose $E_1$ $\k$-orthogonal to~$E_2$ and~$E_3,$ with $\k(E_1,E_1)=1,$ which is possible because of
the signature of $\k.$ Hence we obtain a $\k$-hyperbolic basis
$(E_1,E_2,E_3).$ Since $u$ is $\k$-symmetric, the matrix of $u$  in
the basis constructed above is
\[
\left(\begin{array}{ccc}
 \gamma &a&0\\
0&1&0\\
a&b &1
\end{array}\right),
\]
whose characteristic polynomial is $p_u(x)=(x-\gamma)(x-1)^2.$ Since
$u$ has only one eigendirection, $\gamma=1$ and the matrix of
$u^{-1}$ is
\[
\left(
\begin{array}{ccc}
1&-a&0\\
0&1&0\\
-a&a^2-b &1
\end{array}\right).
\]
The Lax pair is given by:
\begin{gather*}
\dot{x}_1 =-ax_1x_2+(a^2-b)x_2^2+(b-1)x_2x_3,\\
\dot{x}_2=a x_2^2,\\
\dot{x}_3=ax_1^2-(a^2-b)x_1x_2-ax_2x_3.
\end{gather*}
Notice that the hypothesis imply that $a\ne 0.$ An idempotent is
produced by setting  $x_2=1/a$ and solving the other two equations,
with this value for $x_2.$ \end{proof}

\begin{lemma}
Let $u$ be diagonalizable. Then the pseudo-Riemannian metric on
${\mathrm{Sl}}(2,\r)$ defined by $u$ has a definite quadratic first
integral (and hence is complete) provided that the Lax pair vector
field has no idempotents.
\end{lemma}

\begin{proof} If $u$ is diagonalizable then there is a
$k$-orthonormal basis $E_1$, $E_2$, $E_3$ with $E_1$, $E_2$, $E_3$ eigenvectors
of~$u.$ In fact, if $u$ has one eigenvalue with eigenspace of
dimension~3, then $u$ is a~multiple of the identity and the
assertion follows trivially. If~$u$ has three dif\/ferent eigenvalues,
then the proof is straightforward. Assume $u$ has two eigenvalues
and that one of them has eigenspace of dimension~2, i.e.\ there
exists a basis $\{E_1,E_2,E_3\}$ such that $u(E_1)=\lambda E_1,$
$u(E_2)=\lambda E_2,$ $u(E_3)=\mu E_3,$ $\lambda\ne\mu.$ Then the
restriction of $u$ to the eigenspace corresponding to $\lambda$ is a
multiple of the identity. Since $\k(E_1,E_3)=\k(E_2,E_3)=0,$ then
$\k(E_3,E_3)\ne 0.$ Hence the $\k$-non degenerate subspace
${\mathrm{Vect}}\{E_3\}$ is invariant by $u.$ This implies that the
eigenspace ${\mathrm{Vect}}\{E_1,E_2\}={\mathrm{Vect}}\{E_3\}^\perp$
is $\k$-non degenerate. And the assertion follows.

Hence the matrix of $u^{-1}$ in the $\k$-orthonormal basis of
eigenvectors is given by
\[
\left(\begin{array}{ccc}
\lambda_1  & 0   & 0 \\
0 &  \lambda_2 & 0\\
0  & 0   &  \lambda_3
\end{array}\right),
\]
where $\lambda_1^{-1}$, $\lambda_2^{-1}$ and $\lambda_3^{-1}$ are the
eigenvalues of $u.$ The Lax pair is given by
%\label{eq:esp1}
\begin{gather*}
\dot{x}_1 = (\lambda_3 -\lambda_2)\,x_2x_3 =ax_2x_3, \\
\dot{x}_2=(\lambda_1 - \lambda_3)\,x_1x_3 = bx_1x_3,\\
\dot{x}_3=(\lambda_1 -\lambda_2)\,x_1x_2  =(a+b)x_1x_2,
\end{gather*}
where $a=\lambda_3 -\lambda_2,$ $b=\lambda_1-\lambda_3.$  The
equation for idempotents is:
\begin{gather*}
x_1=a\,x_2x_3, \qquad
x_2= b\,x_1x_3,\qquad
x_3=(a+b)\,x_1x_2,
\end{gather*}
which is equivalent to:
\begin{equation}\label{eq:idempotentes}
b(a+b)\,x_1^2=a(a+b)\,x_2^2=ab\,x_3^2.
\end{equation}

Equation (\ref{eq:idempotentes}) has a non trivial solution if and
only if $a,$ $b$ and $a+b$ have the same signature, that is, if and
only if $ab>0.$ Hence the Euler product has no non trivial
idempotents if and only if $ab<0$ or $ab(a+b)= 0.$ If $ab(a+b)=0$
the Lax pair is linearizable. Assume  $ab<0.$ Clearly
\[
bx_1^2-ax_2^2
\]
is a f\/irst integral. This quadratic, which is a semi-def\/inite
quadratic form,  and the quadratic f\/irst integral
$x_1^2+x_2^2-x_3^2$ can be combined into a quadratic positive
def\/inite f\/irst integral.
\end{proof}

\begin{lemma}\label{lm:3.3}
If $u$ has two linearly independent eigenvectors corresponding to
two different eigenvalues then the left invariant pseudo-Riemannian
metric on $\mathrm{Sl}(2,\r)$ defined by $u$ is complete whenever
the Lax pair vector field has no idempotents.
\end{lemma}

\begin{proof} Let $p_u(x)=(x-\alpha)^2(x-\beta)$ be the
characteristic polynomial of $u,$ $E$ an eigenvector of $u$
corresponding to the simple eigenvalue $\beta.$ Denote by $V$ the
subspace $\{ x\in\g : (u-\alpha\,{\mathrm{Id}})^2x=0 \}.$ Then $V$
and the subspace generated by $E$ are  $u$ invariant and $E\perp_k
V.$ The last assertion follows from the following calculation: when
$x\in V,$
\[
u^2(x)-2\alpha u(x)+\alpha^2x=0
\]
and
\begin{gather*}
\alpha^2\k(E,x) = -\k(E,u^2(x)-2\alpha u(x)) = - \k(E,u^2(x))+2\alpha \k(E,u(x))\\
\hphantom{\alpha^2\k(E,x)}{} = -\beta^2 k(E,x) + 2\alpha\beta k(E,x)\
\end{gather*}
since $u$ is $\k$-symmetric. Hence
\[
(\alpha^2 +\beta^2  - 2\alpha\beta )\k(E,x)=0
\]
and $\k(E,x)=0$. (By hypothesis, $\alpha\ne\beta $.)

For every $E'\in V$ eigenvector corresponding to $\alpha $,
$\k(E',E')=0$ because if $\k(E',E')\ne 0,$ then
${\rm{Vect}}\{E,E'\}$ is a $u$ invariant, $\k$~-- non degenerate
subspace. Hence ${\rm{Vect}}\{E,E'\}^{\perp_\k}$ is $u$ invariant of
dimension 1. This implies the existence of an eigenvector $E''\notin
{\rm{Vect}}\{E,E'\}$  which by hypothesis is not possible.

The subspace $V$ is $\k$-non degenerate, of dimension 2 and has an
isotropic vector. Hence it is a hyperbolic plane. Since $\k$ is of
index  $(2,1),$ $\k(E,E)>0.$ Let $E_1=E/\sqrt{\k(E,E)}.$

Let $E_2, E_3\in V$ such that
\begin{gather*}
\k(E_2,E_2)=0=\k(E_3,E_3),\qquad \k(E_2,E_3)=1 \qquad {\mbox{and}}
\\
u(E_2)=\alpha E_2 , \qquad u(E_3)=\gamma E_2 +\alpha E_3,
\end{gather*}
where $0\ne \gamma=\langle E_3, E_3\rangle.$ In this $\k$-hyperbolic
basis, the matrix of $u^{-1}$ is{\samepage
\[
\left(\begin{array}{ccc}
b&0&0\\
0&a & -\gamma a^2 \\
0&0&a
\end{array}\right),
\]
where $a=\alpha^{-1}$ and $b=\beta^{-1}.$}

The Lax pair associated to the pseudo-Riemannian metric is
\begin{gather*}
\dot{x}_1 = \gamma a^2x_3^2, \\
\dot{x}_2= -\gamma a^2\,x_1x_3 +(a-b)\,x_1x_2,\\
\dot{x}_3 = (b-a)\,x_1x_3.
\end{gather*}

The f\/irst and third equations conform a planar system which has
idempotents if and only if
\[
\gamma(b-a) > 0.
\]
The idempotent obtained produces an idempotent for the Lax pair
vector f\/ield. Then, if the Lax pair vector f\/ield has no idempotent
we get a semi-def\/inite quadratic f\/irst integral for the whole system
and a def\/inite quadratic f\/irst integral for the planar system:
\[
(a-b)x_1^2+\gamma a^2x_3^2.
\]
Thus the planar system is complete. The second equation is of the
type
\[
\dot{x}_2= f(t)+ g(t)x_2.
\]
where $f$, $g$ are def\/ined for all $t,$ therefore it is also complete.
\end{proof}

\begin{lemma}
If $u$ has two linearly independent eigenvectors corresponding to
the same eigenvalue, then the Lax pair associated to $u$ on
$\mathrm{Sl}(2,\r)$ is linearizable.
\end{lemma}

\begin{proof} Assuming that $u$ has two eigendirections and exactly
one eigenvalue, there is a basis of $\g$, $E_1$, $E_2$, $E_3$, such that
\[
M_{{\cal{B}}}(u)=\left(\begin{array}{ccc}
\alpha  & 1   & 0 \\
0 &  \alpha  & 0\\
0  & 0   &  \alpha
\end{array}\right).
\]
We have that $E_1^\perp={\mathrm{Vect}}\{E_1,E_3\}.$ In fact, since
$u$ is $\k$-symmetric,
\[
\begin{array}{ccc}
\k(E_1,u(E_2))&=&\k(u(E_1),E_2)\vspace{-1mm}\\
\shortparallel&&\shortparallel\vspace{-1mm}\\
\k(E_1,E_1)+ a\k(E_1,E_2)&=& a\k(E_1,E_2)
\end{array}
\]
which implies that $\k(E_1,E_1)=0,$ and
\[
\begin{array}{ccc}
\k(E_3,u(E_2))&=&\k(u(E_3),E_2)\vspace{-1mm}\\
\shortparallel&&\shortparallel\vspace{-1mm}\\
\k(E_3,E_1)+ a\k(E_3,E_2)&=& a\k(E_3,E_2)
\end{array}
\]
which implies that $\k(E_1,E_3)=0.$

Hence
\[
\k(E_1,E_2)\ne 0.
\]
The subspace $V={\rm{Vect}}\{E_1,E_2\}$ is non degenerate and $u$
invariant. Thus $V^{\perp_\k}$ is a one-dimensional subspace, non
degenerate and $u$ invariant.

The Lax pair associated to the pseudo-Riemannian metric is
\begin{gather*}
\dot{x}_1=a^2x_3^2, \qquad
\dot{x}_2=-a^2\,x_1x_3,\qquad
\dot{x}_3=0,
\end{gather*}
which is linearizable.
\end{proof}

We can summarize the situation in $\mathrm{Sl}(2,\r)$ with the
following proposition

\begin{proposition}\label{pr:sl}
A left invariant pseudo-Riemannian metric in $\mathrm{Sl}(2,\r)$
with associated isomorphism $u$ and minimal polynomial $q_u(x)$, is
complete if and only if
\begin{itemize}\itemsep=0pt
\item[$(i)$] $q_u(x)$ has degree less than or equal to $2$. In this case the
Euler equation is linearizable;
\item[$(ii)$] $q_u(x)=(x-\alpha_1)(x-\alpha_2)(x-\alpha_3)$ where $\alpha_i\ne \alpha_j$ if $i\ne j$
and
\[
\left(\frac{1}{\alpha_3}-\frac{1}{\alpha_2}\right)
\left(\frac{1}{\alpha_3}-\frac{1}{\alpha_1}\right)> 0,
\]
where $\alpha_3$ is the eigenvalue of $u$ corresponding to an
eigenvector $E$ such that $\k(E,E)<0.$ In this case, the Euler
equation of a complete left invariant pseudo-Riemannian metric has a
definite quadratic first integral;
\item[$(iii)$] $q_u(x)=(x-\alpha_1)(x-\alpha_2)^2$ and
\[
\langle
E_3,E_3\rangle\left(\frac{1}{\alpha_2}-\frac{1}{\alpha_1}\right)
>0,
\]
where $E_3$ is a cyclic vector associated to $\alpha_2.$
\end{itemize}
\end{proposition}

\section{Some examples and remarks}\label{section3}

It is worth noting that, except for the case ${\mathrm{E}}(2)$ and
the case in Proposition~\ref{pr:sl}$(iii)$ we have shown that a left
invariant pseudo-Riemannian metric on a non unimodular Lie group is
complete if and only if either the Euler equation is linearizable or
has a def\/inite quadratic f\/irst integral. The following example shows
that this is not the case for ${\mathrm{E}}(2)$ (and an example can
be found in the same way in ${\mathrm{Sl}}(2,\r))$. However in these
cases the dif\/ferential equation is of the type:
\begin{gather*}
\dot{x}_1= X_1(x_1,x_2), \qquad
\dot{x}_2= X_2(x_1,x_2),\qquad
\dot{x}_3= X_3(x_1,x_2,x_3)
\end{gather*}
($X_1$, $X_2$, $X_3$ quadratic polynomials) where the two f\/irst equations
have a def\/inite quadratic f\/irst integral.

%\subsection*{Example 1}
\begin{example}\label{example1}
Consider on ${\mathrm{E}}(2)$ the pseudo-Riemannian metric given by
\[
-x_1^2+x_2^2+2x_2x_3.
\]
in the basis considered in Proposition~\ref{proposition2}. The
corresponding Euler equation is
\begin{gather*}
\dot{\xi}_1=\xi_1\xi_2,\qquad
\dot{\xi}_2=\xi_3^2,\qquad
\dot{\xi}_3=-\xi_2\xi_3.
\end{gather*}
This equation is complete and is not linearizable. However $\xi_1$
is not bounded, hence no f\/irst integral can have bounded level
surfaces.

 In \cite{GL} a characterization of the
complete pseudo-Riemannian metrics on~${\mathrm{Sl}}(2,\r)$ is
given:

\medskip

\noindent {\bf Theorem\footnote{\noindent
{\bf{Theorem.}} {\it Si la m\'etrique invariante \`a gauche sur
${\mathrm{Sl}}(2,\r)$ d\'efinie par l'endomorphisme $\phi$ est
lorentzienne, elle n'est compl\`ete que dans les deux cas suivants:
\begin{enumerate}\itemsep=0pt
\item[$1)$] $\phi$ admet un sous-espace propre de dimension $2$ au moins;
\item[$2)$] $\phi$ est diagonalizable sur $\r,$ admet deux valeurs propres
oppose\'ees, et les vecteurs propres associ\'es engendrent un plan
lorentzien.
\end{enumerate}}}.} {\it The Lorentzian left invariant pseudo-Riemannian metric on
${\mathrm{Sl}}(2,\r)$ is defined by the isomorphism $\phi$ is
complete only in the two following cases:
\begin{enumerate}\itemsep=0pt
\item[$1)$] $\phi$ has an eigenspace of dimension at least two;
\item[$2)$] $\phi$ is diagonalizable over the real numbers, has two eigenvalues of different sign
and the eigenvectors associated to these eigenvalues generate a
Lorentzian plane.
\end{enumerate}}

However, the following examples show that this characterization is
not correct.
\end{example}

%\subsection*{Example 2}
\begin{example}\label{example2}
Consider on ${\mathrm{Sl}}(2,\r)$ the left invariant Lorentzian
metric given in a $\k$-orthonormal basis ${\cal{B}}=(E_1,E_2,E_3)$
by
\[
x_1^2 + x_2^2+x_2x_3.
\]
In the $\k$-hyperbolic basis deduced from ${\cal{B}},$
${\cal{B}}^\prime=(E_1,
(1/\sqrt{2})(E_2+E_3),(1/\sqrt{2})(E_2-E_3)),$ the $\k$-symmetric
isomorphism $u$ that corresponds to the Lorentzian metric is given
by the matrix
\[
M_{{\cal{B}}^\prime}(u)=\left( \begin{array}{ccc} 1&0&0\\
                                           0&1/2&1\\
                                           0&0&1/2
                        \end{array}\right).
\] The characteristic polynomial of $u$ is $p_u(x)=(x-1/2)^2(x-1),$
with eigendirections  $E_1$ and $E_2+E_3.$ The Lax pair equation
corresponding to this Lorentzian metric in ${\cal{B}}^\prime$
coordinates is
\begin{gather*}
\dot{x_1}=4x_3^2,\qquad
\dot{x_2}=x_1(x_2-4x_3),\qquad
\dot{x_3}=-x_1x_3.
\end{gather*}
The planar system given by the  f\/irst and third  equations is
complete. Hence the system is complete. Notice that $u$ has neither
an eigenspace of dimension 2 nor is it diagonalizable.
\end{example}

%\subsection*{Example 3}
\begin{example}\label{example3}
Let $u$ be the isomorphism with matrix
\[
M_{{\cal{B}}}(u)=\left( \begin{array}{ccc} 1/2&0&0\\
                                           0&1/3&0\\
                                           0&0&1
                        \end{array}\right).
\]
The pseudo-Riemannian metric is given by:
\[
\langle x,x\rangle=\k(u(x),x)=\tfrac{1}{2}x_1^2+\tfrac{1}{3}x_2^2-x_3^2
\]
which is Lorentzian. Notice that all eigenvalues of $u$ are positive

The Lax pair corresponding to this pseudo-Riemannian metric is
\begin{gather*}
\dot{x}_1= -2x_2x_3, \qquad
\dot{x}_2=x_1x_3,\qquad
\dot{x}_3=-x_1x_2.
\end{gather*}
This equation has the following f\/irst integrals
\[
x_1^2+x_2^2-x_3^2\qquad  \mbox{and} \qquad x_2^2+x_3^2.
\]
Hence the solutions are bounded and the equation is complete.
\end{example}

%\subsection*{Example 4}
\begin{example}\label{example4}
Consider the pseudo-Riemannian metric in the same basis as
before given by
\[
x_1^2-x_2^2+2x_3^2.
\] The isomorphism $u$
associated to this pseudo-Riemannian metric has $E_1$, $E_2,$ and
$E_3$ as eigenvectors with (respectively) eigenvalues $1$, $-1$ and~$2$.
Hence we have two opposite eigenvalues whose corresponding
eigenvectors generate a Lorentzian plane. The Lax pair equation in
this case is:
\begin{gather*}
\dot{x}_1= \tfrac{1}{2}x_2x_3,\qquad
\dot{x}_2=\tfrac{3}{2}x_1x_3,\qquad
\dot{x}_3=2x_1x_2.
\end{gather*}
Since $\big(\sqrt{\tfrac{1}{3}}, 1, 2\sqrt{\tfrac{1}{3}}\big)$ is
an idempotent for the Lax pair vector f\/ield, the pseudo-Riemannian
metric is not complete.

We have shown that for unimodular Lie groups of dimension 3 the
simplest case for incompleteness occur. This is not the case in
general as the following example shows.
\end{example}

%\subsection*{Example 5}
\begin{example}\label{example5}
Let $\g$ be a non unimodular Lie algebra de Lie of dimension 3. Then
there is a~basis $(e_1,e_2,e_3)$, of $\g$ such that the Lie bracket
is given by
\[
\begin{array}{c|c|c|c}
\rule[-2mm]{0cm}{4mm}[\, ,\,]  & e_1  & e_2   & e_3 \\ \hline
\rule[-2mm]{0cm}{6mm} e_1
&  0   & \alpha e_2+\beta e_3  & \gamma e_2 +\delta e_3\\
\hline \rule[-2mm]{0cm}{6mm} e_2
&\phantom{\alpha e_2+\beta e_3}  &      0               & 0\\
\hline \rule[-2mm]{0cm}{6mm} e_3      &      & &         0
\end{array}
\]
where $\alpha+\delta = 2$~\cite{M}. For $\alpha=1/2$, $\delta
=3/2$, and energy $2\xi_1^2-\xi_2\xi_3$, the Euler equation is:
\begin{gather*}
\dot{\xi}_1=\xi_2\,\xi_3,\qquad
\dot{\xi}_2= \xi_1\xi_2,\qquad
\dot{\xi}_3=3 \xi_1\xi_3.
\end{gather*}

Clearly the Euler vector f\/ield has no idempotents and the curve
\[
\xi(t)=\left(\dfrac{1}{\sqrt{2}}\dfrac{1}{1-\sqrt{2}t},\dfrac{1}{(1-\sqrt{2}t)^{1/2}},
\dfrac{1}{(1-\sqrt{2}t)^{3/2}}\right)
\]
is an incomplete solution with initial condition
$\xi(0)=(1/\sqrt{2},1,1)$.
\end{example}

\subsection*{Acknowledgements}

The authors wish to thank the anonymous referee for the careful
reading and the pertinent observations that led, in particular, to a
revision of the proof of Lemma~\ref{lm:3.3} and to reformulate
accordingly Proposition~\ref{pr:sl}.

\pdfbookmark[1]{References}{ref}
\LastPageEnding

\end{document}